\documentclass[12pt] {article}
\usepackage{amsmath,amsthm}
\usepackage{amssymb,latexsym}
\usepackage{graphicx}
\usepackage{amssymb}
\usepackage{amsfonts}
\usepackage{amscd}

\usepackage{amscd}
\title{Theory of valuations on manifolds: a survey.}
\date{}
\author{ Semyon Alesker \footnote{Partially supported by ISF grant 1369/04.}
\\  { \normalsize Department of Mathematics, Tel Aviv University, Ramat Aviv}
 \\  { \normalsize 69978 Tel Aviv,
Israel }
\\ {\normalsize e-mail: semyon@post.tau.ac.il}}
\def\One{{1\hskip-2.5pt{\rm l}}}
\def\eps{\varepsilon}
\def\alp{\alpha}
\def\ome{\omega}

\def\lam{\lambda}

\def\to{\rightarrow}

\def\agrc{{}\!^ {\textbf{C}} {\cal A}Gr}

\def\pt{\partial}

\def\RR{\mathbb{R}}
\def\CC{\mathbb{C}}

\def\ZZ{\mathbb{Z}}

\swapnumbers
\newtheorem{theorem}{Theorem}[section]

\newtheorem{examples}[theorem]{Examples}

\newtheorem{proposition}[theorem]{Proposition}

\theoremstyle{definition}

\newtheorem{definition}[theorem]{Definition}
\newtheorem{remark}[theorem]{Remark}
\theoremstyle{conjecture}
\newtheorem{conjecture}[theorem]{Conjecture}
\theoremstyle{proposition-definition}
\newtheorem{proposition-definition}[theorem]{Proposition-Definition}

\def\cf{{\cal F}}

  \def\cf{{\cal F}}
  
 \def\ck{{\cal K}} 
  
\def\cp{{\cal P}}

\def\cv{{\cal V}} \def\cw{{\cal W}} 
 
\def\inj{\hookrightarrow}

\newcommand \supp{\operatorname{supp} \,}

\newcommand \codim{\operatorname{codim} \,}
\begin{document}
\maketitle \setcounter{section}{-1}
\begin{abstract}
This is a non-technical survey of a recent theory of valuations on
manifolds constructed in \cite{part1}-\cite{part3} and actually a
guide to this series of articles. We review also some recent
related results obtained by a number of people. We formulate some
open questions.
\end{abstract}
\section{Introduction.} In convexity there are many
geometrically interesting and well known examples of valuations on
convex sets: Lebesgue measure, the Euler characteristic, the
surface area, mixed volumes, the affine surface area. For a
description of older classical developments on this subject we
refer to the surveys McMullen-Schneider \cite{mcmullen-schneider},
McMullen \cite{mcmullen-survey}. For the general background on
convexity we refer to the book by Schneider \cite{schneider-book}.

Approximately during the last decade there was a significant
progress in this classical subject which has led to new
classification results of various classes of valuations, to
discovery of new structures on them. This progress has shed a new
light on the notion of valuation which allowed to generalize it in
some cases to the more general setting of valuations on manifolds
and on not necessarily convex sets (a concept which in any case
has no meaning on a general manifold). On the other hand, the
author's feeling is that the notion of valuation equips smooth
manifolds with  a new general rich structure. Valuations on
manifolds were introduced and studied in a series of four
articles: \cite{part1}, \cite{part2}, \cite{part4} by the author,
and \cite{part3} by J. Fu and the author. This theory depends
heavily on and is a continuation of the classical theory of
valuations on convex subsets of an affine space; it combines
probably most of the results obtained on translation invariant
continuous valuations on convex sets. Also the notions of normal
(or characteristic) cycle and Legendrian currents and tools from
geometric measure theory to study them (see e.g.
\cite{fu-89}-\cite{fu-94}), turned out to be very useful in this
new theory, as well as tools from representation theory.


The goal of this article is to give a non-technical overview of
\cite{part1}, \cite{part2}, \cite{part3}, \cite{part4}. We also
mention a few recent closely related results by a number of
people. We state a number of open questions on valuations on
manifolds.

In \cite{part2} the notion of smooth valuation on a smooth
manifold was introduced. Roughly put, a smooth valuation can be
thought as a finitely additive $\CC$-valued measure on a class of
nice subsets; this measure is required to satisfy some additional
assumptions of continuity (or rather smoothness) in some sense.
The basic examples of smooth valuations on a general manifold $X$
are smooth measures on $X$ and the Euler characteristic. Moreover,
the well known intrinsic volumes of convex sets can be generalized
to provide examples of smooth valuations on an arbitrary
Riemannian manifold; these valuations are known as
Lipschitz-Killing curvatures.

The article is organized as follows. In Section \ref{overview} we
give an overview of some necessary previously known facts on
translation invariant valuations on convex sets.  In Section
\ref{more} we have surveyed relevant results on translation
invariant valuations on convex sets due to the author but which
were proved in other places than \cite{part1}-\cite{part3}, and
due to Bernig-Br\"ocker \cite{bernig-brocker}, and J. Fu
\cite{fu-unitary}. Section \ref{main-results} contains the
description of \cite{part1}-\cite{part3}; this is the main section
of this article. Section \ref{guide} plays a role of an appendix:
this is a short guide to \cite{part1}-\cite{part3} where we
indicate in what part of the series \cite{part1}-\cite{part3} one
can find the proofs of the results discussed in Section
\ref{main-results}.

This article is not a survey of the developments of valuation
theory, not even during the last decade, and there are a number of
interesting developments which are not discussed here. To mention
just a few of them, these are: \cite{alesker-geom-dedic},
\cite{gates-hug-schneider-05}, \cite{klain-95}, \cite{klain-2000},
\cite{ludwig-99}, \cite{ludwig-03}, \cite{ludwig-reitzner-99},
\cite{schneider-96}. Some of them (particularly Klain
\cite{klain-95}, Schneider \cite{schneider-96}) were very
influential on the progress discussed in this article.

{\bf Acknowledgements.} I express my gratitude to J. Bernstein and
V.D. Milman for numerous useful discussions and their attention to
my work over the years. I thank J. Fu for his collaboration with
me. Also I thank the referee for numerous remarks on the first
version of the article.

\section{A brief overview of valuations on convex
sets.}\label{overview} In this section we present a brief overview
of necessary facts from the classical theory of valuations on
convex sets. Let $V$ be a finite dimensional real vector space of
dimension $n$. Let us denote by $\ck(V)$ the family of non-empty
convex compact subsets of $V$. Then $\ck(V)$ has a natural
topology. To define it, let us fix a Euclidean metric on $V$. The
Hausdorff metric $d_H$ on $\ck(V)$ is defined as follows:
$$d_H(A,B):=\inf\{\eps>0|\,A\subset (B)_\eps \mbox{ and } B\subset
(A)_\eps\}$$ where $(U)_\eps$ denotes the $\eps$-neighborhood of a
set $U$. It is well known (the Blaschke selection theorem) that
$\ck(V)$ equipped with the Hausdorff metric $d_H$ is a locally
compact space. If we choose a different Euclidean metric on $V$
the corresponding Hausdorff metric will define the same topology
on $\ck(V)$.
\begin{definition}\label{defval}
A scalar valued functional
$$\phi\colon \ck(V)\to \CC$$
is called a {\itshape convex valuation} if
$$\phi(A\cup B)=\phi(A)+\phi(B)-\phi(A\cap B)$$
whenever $A,B,A\cup B\in \ck(V)$.
\end{definition}
\begin{remark}
In Definition \ref{defval} the notion we called {\itshape convex
valuation} is called just {\itshape valuation} in all the
classical literature. We have made this change of terminology in
order to emphasize on one hand that in the sequel we will consider
valuations defined on not necessarily convex sets, and on the
other hand that the new notion of valuation discussed below
generalizes in a sense the classical notion from Definition
\ref{defval}.
\end{remark}
\begin{definition}
A convex valuation $\phi$ is called {\itshape continuous} if
$\phi$ is continuous in the Hausdorff metric.
\end{definition}

Let us denote by $Val(V)$ the space of {\itshape translation
invariant} continuous convex valuations. Equipped with the
topology of uniform convergence on compact subsets of $\ck(V)$ the
space $Val(V)$ becomes a Banach space. Let us give some examples:

1) A Lebesgue measure belongs to $Val(V)$.

2) The Euler characteristic $\chi$ belongs to $Val(V)$ (recall
that $\chi(K)=1$ for any $K\in \ck(V)$).

3) Fix a convex compact set $A\in \ck(V)$. Let $vol$ denotes a
Lebesgue measure on $V$. Define
$$\phi(K):= vol(K+A)$$
where $K+A$ denotes the Minkowski sum of $K$ and $A$, namely
$K+A:=\{k+a|\, k\in K,a\in A\}$. Then $\phi\in Val(V)$.

4) Let us fix an integer $i=0,1,\dots,n$. Let us fix
$A_1,\dots,A_{n-i}\in \ck(V)$. Then the mixed volume $K\mapsto
V(K[i],A_1,\dots,A_{n-i})$ belongs to $Val(V)$, where $K[i]$ means
that $K$ is taken $i$ times (for the notion of mixed volume see
e.g. Schneider's book \cite{schneider-book}).

5) There is a different construction of continuous translation
invariant convex valuations based on the theory of complex and
quaternionic plurisubharmonic functions; see
\cite{alesker-jdg-03}, \cite{alesker-adv-05}.

It was conjectured by P. McMullen \cite{mcmullen-80} that linear
combinations of mixed volumes are dense in $Val(V)$. This
conjecture was proved by the author in \cite{alesker-gafa-01} in a
much stronger form. In order to formulate this result, called
Irreducibility Theorem, let us remind a few necessary facts.

Let $\phi$ be a convex valuation. Let $\alp$ be a complex number.
We say that $\phi$ is $\alp$-homogeneous if
$$\phi(\lam K)=\lam^\alp\phi(K) \mbox{ for any } \lam>0,\, K\in
\ck(V).$$ Let us denote by $Val_\alp(V)$ the subspace of $Val(V)$
of $\alp$-homogeneous convex valuations. For instance, the Euler
characteristic is 0-homogeneous, a Lebesgue measure is
$n$-homogeneous, and the mixed volumes from Example 4) are
$i$-homogeneous. The following result is due to P. McMullen
\cite{mcmullen-euler}.
\begin{theorem}[\cite{mcmullen-euler}]
Let $n=\dim V$. Then
$$Val(V)=\bigoplus_{i=0}^n Val_i(V).$$
\end{theorem}

\begin{remark}
1) It is easy to see that $Val_0(V)$ is one dimensional and is
spanned by the Euler characteristic $\chi$.

2) It was shown by Hadwiger \cite{hadwiger-book} that $Val_n(V)$
is also one dimensional and is spanned by a Lebesgue measure.
\end{remark}

We say that a convex valuation $\phi$ is {\itshape even} if for
any $K\in\ck(V)$ one has $\phi(-K)=\phi(K)$. Similarly $\phi$ is
called {\itshape odd} if for any $K\in\ck(V)$ one has
$\phi(-K)=-\phi(K)$. The space of $i$-homogeneous valuations is
decomposed further into a direct sum with respect to parity:
$$Val_i(V)=Val_i^{ev}(V)\oplus Val_i^{odd}(V)$$
where the notation is obvious.

The group $GL(V)$ of invertible linear transformations of $V$ acts
naturally on $Val(V)$ as follows:
$$(g\phi)(K)=\phi(g^{-1}(K)) \mbox{ for } g\in GL(V),\phi\in Val(V), K\in
\ck(V).$$ This action is continuous and preserves the degree of
homogeneity and the parity of valuations. The Irreducibility
Theorem says the following.
\begin{theorem}[\cite{alesker-gafa-01}]\label{irrthm}
The natural representation of $GL(V)$ in $Val_i^{ev}(V)$ and
$Val_i^{odd}(V)$ is irreducible for any $i=0,1,\dots,n=\dim V$
(i.e. there is no proper closed $GL(V)$-invariant subspace).
\end{theorem}
This result implies P. McMullen's conjecture (see
\cite{alesker-gafa-01}). The proof used most of the known (by that
time) results on translation invariant continuous convex
valuations in combination with tools from representation theory
\cite{beilinson-bernstein}.

The following definition is a special case of a more general
notion of a smooth vector in a representation space of a group.
\begin{definition}\label{smoothdef}
A convex valuation $\phi\in Val(V)$ is called {\itshape smooth} if
the map $GL(V)\to Val(V)$ defined by $g\mapsto g(\phi)$ is
infinitely differentiable.
\end{definition}
Let us denote by $Val^\infty(V)$ the subset of smooth convex
valuations in the sense of Definition \ref{smoothdef}. It is well
known (for general representation theoretical reasons) that
$Val^\infty(V)\subset Val(V)$ is a linear $GL(V)$-invariant
subspace dense in $Val(V)$. The space $Val^\infty(V)$ has a
canonical Fr\'echet topology. Denoting by $Val_i^\infty(V)$,
$Val_i^{ev,\infty}(V)$, $Val_i^{odd,\infty}(V)$ the subspaces of
smooth vectors in $Val_i(V)$, $Val_i^{ev}(V)$,$Val_i^{odd}(V)$,
respectively, one easily deduces the following versions of
McMullen's decomposition:
\begin{eqnarray}\label{mcdecomp1}
Val^\infty(V)=\bigoplus_{i=0}^nVal_i^\infty(V);\\\label{mcdecomp2}
Val_i^\infty(V)=Val_i^{ev,\infty}(V)\oplus Val_i^{odd,\infty}(V).
\end{eqnarray}
\section{Some results and questions on translation invariant valuations on convex sets.}\label{more}
In this section we review a number of results on translation
invariant (convex) valuations due to several people which are very
closely related to the material discussed in Section
\ref{main-results} below. For an $n$-dimensional real vector space
$V$ we denote by $Val(V)$ the space of continuous translation
invariant {\itshape convex} valuations in the sense of Definition
\ref{defval}, and by $Val^\infty(V)$ the subspace of smooth convex
valuations in the sense of Definition \ref{smoothdef}. It turns
out that valuations from $Val^\infty(V)$ have the following
important property (compare with Proposition \ref{comparison}
below): they can be naturally evaluated on compact sets of
{\itshape positive reach} (or, in other terminology, semi-convex
sets); this class of sets contains all the compact convex sets and
all compact submanifolds with corners.


Recall that we have McMullen's decomposition (\ref{mcdecomp1})
with respect to the degree of homogeneity:
\begin{eqnarray}\label{grading}
Val^\infty(V)=\bigoplus_{i=0}^n Val_i^\infty(V).
\end{eqnarray}
The following result was proved in \cite{alesker-gafa-04}.
\begin{theorem}\label{poincare}
The space $Val^\infty(V)$ has a canonical continuous product
$Val^\infty(V)\times Val^\infty(V)\to Val^\infty(V)$. Then
$Val^\infty(V)$ becomes a graded algebra:
$$Val_i^\infty(V)\cdot Val_j^\infty(V)\subset Val_{i+j}^\infty(V).$$
Moreover, it satisfies the Poincar\'e duality, namely for any
$i=0,1,\dots,n$ the product
$$Val_i^\infty(V)\times Val_{n-i}^\infty(V)\to Val_n^\infty(V)(\, =\CC\cdot vol)$$
is a perfect pairing, in other words, the induced map
$$Val_i^\infty(V)\to (Val_{n-i}^\infty(V))^*\otimes
Val_n^\infty(V) $$ is injective and has a dense image in the weak
topology.
\end{theorem}
Thus $Val^\infty(V)$ is a graded algebra satisfying the Poincar\'e
duality; such algebras are often called Frobenius algebras. The
proof of this theorem, besides the construction of the product
from \cite{alesker-gafa-04}, uses the full generality of the
Irreducibility Theorem.

The following result is a version of the hard Lefschetz theorem
for even valuations. In order to formulate it, let us fix a
Euclidean metric on $V$. Let us denote by $V_1\in Val_1^\infty(V)$
the first intrinsic volume on $V$ (see e.g. \cite{schneider-book},
p. 210); by Hadwiger's theorem \cite{hadwiger-book} this is, up to
a constant, the only non-zero 1-homogeneous $SO(n)$-invariant
translation invariant continuous valuation (which is automatically
smooth, see Proposition \ref{prop-smooth} below).
\begin{theorem}[\cite{alesker-gsn-04}]\label{hlmult}
Fix an integer $i$, $0\leq i<n/2$. Then the operator
$$Val_i^{ev,\infty}(V)\to Val_{n-i}^{ev,\infty}$$
defined by $\phi\mapsto V_1^{n-2i}\cdot \phi$ is an isomorphism.
\end{theorem}
The proof of this theorem in \cite{alesker-gsn-04} uses, besides
much of the machinery of the valuation theory, the results on the
Radon transform on Grassmannians due to Gelfand, Graev, Ro\c su
\cite{gelfand-graev-rosu} and the solution of the cosine transform
problem for Grassmannians by Bernstein and the author
\cite{alesker-bernstein}. The analogous result was conjectured in
\cite{alesker-gsn-04} for odd valuations:
\begin{conjecture}[\cite{alesker-gsn-04}]\label{conjecture}
Fix an integer $i$, $0\leq i<n/2$. Then the operator
$$Val_i^{odd,\infty}(V)\to Val_{n-i}^{odd,\infty}$$
defined by $\phi\mapsto V_1^{n-2i}\cdot \phi$ is an isomorphism.
\end{conjecture}
We would like to state another version of the hard Lefschetz
theorem for valuations which was proved by the the author in
\cite{alesker-jdg-03} in the even case, and by Bernig-Br\"ocker
\cite{bernig-brocker} in full generality very recently. Let us
consider the operator
$$L\colon Val^\infty(V)\to Val^\infty(V)$$
defined by
\begin{eqnarray}\label{L}
(L\phi)(K):=\frac{d}{d\eps}|_{\eps=0}\phi(K+\eps\cdot D)
\end{eqnarray}
where $D$ is the unit Euclidean ball. (Note that by a result of P.
McMullen \cite{mcmullen-euler}, $\phi(K+\eps\cdot D)$ is a
polynomial in $\eps\geq 0$ of degree at most $n$.) It is easy to
see that the operator $L$ decreases the degree of homogeneity by
one. Then one has
\begin{theorem}\label{hlt-notm}
Fix an integer $i$, $n/2< i\leq n$. Then
$$L^{2i-n}\colon Val_{i}^\infty(V)\to Val^\infty_{n-i}(V)$$
is an isomorphism.
\end{theorem}
Note that the author's proof \cite{alesker-gsn-04} in the even
case used the tools from integral geometry (such as Radon and
cosine transforms on Grassmannians), while Bernig and Br\"ocker
\cite{bernig-brocker} used their description of the forms on the
cotangent bundle defining the zero valuation (see Remark
\ref{rmk-kernel} (2) of this article) in terms of the Rumin
operator, in combination with notion of Laplacian on translation
invariant valuations they introduced.

\begin{remark}\label{polef}
The terminology "Poincar\'e duality" and "hard Lefschetz theorem"
comes from the formal analogy of these properties of the algebra
of valuations with the corresponding properties of the cohomology
algebra of compact K\"ahler manifolds. The Poincar\'e duality is
one of the most basic properties of general compact oriented
manifolds, and the hard Lefschetz theorem is one of the most basic
properties of general compact K\"ahler manifolds (see e.g.
\cite{griffiths-harris}).
\end{remark}

Now let us discuss valuations invariant under a group. Let us fix
from now on a Euclidean metric on $V$. Let $G$ be a {\itshape
compact} subgroup of the orthogonal group. Let us denote by
$Val^G(V)$ the subspace of $Val(V)$ of $G$-invariant convex
valuations. One has the following result.
\begin{proposition}
The space $Val^G(V)$ is finite dimensional if and only if $G$ acts
transitively on the unit sphere of $V$.
\end{proposition}
The "if" part of this proposition was proved in
\cite{alesker-adv-00}, Theorem 8.1. The "only if" part is
announced in print for the first time here; its proof is not very
difficult. Thus in the case when $G$ acts transitively on the unit
sphere one may hope to obtain an explicit finite classification
list of valuations $Val^G(V)$. Thus from now on we will assume
that {\itshape $G$ acts transitively on the unit sphere}.
\begin{proposition}[\cite{alesker-jdg-03}; \cite{alesker-gafa-04},
Theorem 0.9(ii)]\label{prop-smooth} Under the assumption that $G$
acts transitively on the unit sphere one has
$$Val^G(V)\subset Val^\infty(V).$$
\end{proposition}
\begin{remark}\label{G-smooth}
It is important to emphasize that Proposition \ref{prop-smooth}
implies that any $G$-invariant convex valuation from $Val^G(V)$
can be naturally evaluated on compact submanifolds with corners,
and in fact on the larger class of compact sets of positive reach.
\end{remark}
Obviously $Val^G(V)$ is a subalgebra, and one has McMullen's
decomposition with respect to the degree of homogeneity
$$Val^G(V)=\bigoplus_{i=0}^n Val_i^G(V).$$
Define
\begin{eqnarray}\label{hnumber}
h_i:=\dim Val_i^G(V).
\end{eqnarray}
Then $Val_0^G(V)$ is spanned by the Euler characteristic, and
$Val_n^G(V)$ is spanned by a Lebesgue measure.
\begin{theorem}\label{thm-G}
Assume that the group $G$ acts transitively on the unit sphere.

(i) (\cite{alesker-gafa-04}, Theorem 0.9) $Val^G(V)$ is a finite
dimensional graded algebra (with the grading given by the degree
of homogeneity) satisfying the Poincar\'e duality, i.e.
$$Val_i^G(V)\times Val_{n-i}^G(V)\to Val_n^G (\, =\CC\cdot vol)$$
is a perfect pairing. In particular $h_i=h_{n-i}$.

(ii) (\cite{alesker-gafa-04}, Theorem 0.9) Moreover,
$h_1=h_{n-1}=1$. $Val_1^G(V)$ is spanned by the first intrinsic
volume $V_1$; $Val_{n-1}^G(V)$ is spanned by the $(n-1)$-st
intrinsic volume $V_{n-1}$.

(iii) (\cite{alesker-gsn-04}) Assume in addition that $-Id\in G$.
Let $i$ be an integer $0\leq i<n/2$. Then $Val^G(V)$ satisfies a
version of the hard Lefschetz theorem: the operator
$$Val_i^G(V)\to Val_{n-i}^G(V)$$
defined by $\phi\mapsto V_1^{n-2i}\cdot \phi$ is an isomorphism.
\end{theorem}
\begin{conjecture}\label{conjecture1}
In Theorem \ref{thm-G}(iii) the assumption $-Id\in G$ is
unnecessary. More precisely, for any compact group $G$ acting
transitively on the unit sphere the operator
$$Val_i^G(V)\to Val_{n-i}^G(V)$$
defined by $\phi\mapsto V_1^{n-2i}\cdot \phi$ is an isomorphism
for any integer $i$, $0\leq i<n/2$.
\end{conjecture}
Of course, Conjecture \ref{conjecture1} is an immediate
consequence of Conjecture \ref{conjecture}.

Let us state another version of the hard Lefschetz theorem for
$Val^G(V)$ which was proved in \cite{alesker-jdg-03} under the
assumption $-Id\in G$ and in \cite{bernig-brocker} in general (and
which is a consequence of Theorem \ref{hlt-notm}).
\begin{theorem}\label{hlt-nm}
Let $G$ be a group acting transitively on the unit sphere. Let
$L\colon Val^\infty(V)\to Val^\infty(V)$ be the operator defined
by (\ref{L}). Let $i$ be an integer, $n/2<i\leq n$. Then
$$L^{2i-n}\colon Val_i^G(V)\to Val_{n-i}^G(V)$$
is an isomorphism. In particular $h_i\geq h_{i+1}$ for $i\geq
n/2$.
\end{theorem}

Let us consider now concrete examples of compact groups $G$ acting
transitively on the unit sphere. The cases when $G$ is equal
either to the full orthogonal group $O(n)$ or to the special
orthogonal group $SO(n)$ are classical, and there is the following
famous result by Hadwiger \cite{hadwiger-book}.
\begin{theorem}[\cite{hadwiger-book}]\label{so}
$Val^{O(n)}(\RR^n)=Val^{SO(n)}(\RR^n)$, and a basis of this space
is
$$\chi, V_1,V_2,\dots,V_{n-1},vol$$ where $V_i$ denotes the $i$-th
intrinsic volume.
\end{theorem}
Given the Hadwiger theorem \ref{so}, it is not hard to describe
the algebra structure of $Val^{O(n)}(\RR^n)=Val^{SO(n)}(\RR^n)$.
One has
\begin{proposition}[\cite{alesker-gafa-04}, Theorem
2.6]\label{so-alg} The morphism of algebras
$$\CC[x]/(x^{n+1})\to Val^{O(n)}(\RR^n)=Val^{SO(n)}(\RR^n)$$
given by $x\mapsto V_1$ defines an isomorphism of graded algebras.
\end{proposition}

Now let us discuss the other compact groups $G$. It is known in
topology that the condition that $G$ acts transitively on the
sphere, is quite restrictive. In particular, there exists an
explicit classification of {\itshape connected} compact Lie groups
acting transitively on the sphere, due to A. Borel \cite{borel1},
\cite{borel2} and Montgomery-Samelson \cite{montgomery-samelson}.
They have obtained the following list:
\begin{eqnarray}\label{g1}
6  \mbox{ infinite series: }  SO(n), U(n), SU(n), Sp(n),
Sp(n)\cdot Sp(1), Sp(n)\cdot U(1);\\\label{g2}
 3 \mbox{ exceptions: } G_2, Spin(7), Spin(9).
\end{eqnarray}

The next case which has been studied in detail is the group
$G=U(m)$ acting on the standard Hermitian space $\CC^m$. Define
$n:=2m=\dim_\RR \CC^m$. Let us denote by $IU(m)$ the group of
isometries of the Hermitian space $\CC^m$ preserving the complex
structure (then $IU(m)=\CC^m \rtimes U(m)$). Let $\agrc _{j} $
denote the Grassmannian of affine complex subspaces of $\CC^m$ of
complex dimension $j$. Clearly $\agrc _{j} $ is a homogeneous
space of $IU(m)$ and it has a unique (up to a constant)
$IU(m)$-invariant measure (called Haar measure). For every
non-negative integers $p$ and $k$ such that $2p\leq k\leq 2m$  let
us introduce the following valuations:
$$U_{k,p}(K)=\int_{E\in \agrc _{m-p}} V_{k-2p}(K\cap E) \,
dE.$$ Then $U_{k,p}\in Val_k^{U(m)}(\CC^m)$.

\begin{theorem}[\cite{alesker-jdg-03}]\label{u-basis}
The valuations $U_{k,p}$ with $0\leq p \leq \frac{min\{k,
2m-k\}}{2}$ form a basis of the space $Val_k^{U(m)}(\CC^m)$.
\end{theorem}
The proof of this result used, besides the even case of the
Irreducibility Theorem, the even case of the hard Lefschetz
theorem (Theorem \ref{hlt-nm}), and representation theoretical
computations of Howe-Lee \cite{howe-lee}. Some applications of
Theorem \ref{u-basis} to integral geometry of complex spaces can
be found in \cite{alesker-jdg-03}.

The description of the algebra structure of $Val^{U(m)}(\CC^m)$
turned out to be a more difficult problem than for the group
$SO(n)$, and the answer is much more interesting. It was obtained
recently by J. Fu \cite{fu-unitary} in terms of generators and
relations. His result is as follows.
\begin{theorem}[\cite{fu-unitary}]\label{u-alg}
The graded algebra $Val^{U(m)}(\CC^m)$ is isomorphic to the graded
algebra $\CC[s,t]/(f_{m+1},f_{m+2})$ where the generators $s$ and
$t$ have degrees 2 and 1 respectively, and $f_j$ is the degree $j$
component of the power series $log(1+s+t)$.
\end{theorem}
Note that valuations invariant under the other groups from the
list (\ref{g1})-(\ref{g2}) have not been classified, with the only
exception of the group $SU(2)$ acting on $\CC^2\simeq \RR^4$. The
explicit basis of $Val^{SU(2)}(\CC^2)$ in geometric terms was
obtained by the author in \cite{alesker-su2}; the algebra
structure of it has not been computed. Some non-trivial examples
of convex valuations invariant under the quaternionic groups
$Sp(m), Sp(m)\cdot Sp(1)$ were constructed by the author in
\cite{alesker-adv-05} using quaternionic plurisubharmonic
functions.

Another direction which has not been studied in detail is the
description of special classes of valuations on non-affine
manifolds (see Section \ref{main-results} below). For instance one
can prove the following result (compare with Proposition
\ref{so-alg}) which is announced here for the first time; the
details will appear elsewhere. Let $X^n$ denote either the
$n$-dimensional sphere or the $n$-dimensional hyperbolic space
with the standard Riemannian metric.
\begin{proposition}
The algebra of isometry invariant smooth valuations on $X^n$ is
isomorphic to the algebra of truncated polynomials
$\CC[x]/(x^{n+1})$ where a generator is the first
Lipschitz-Killing curvature.
\end{proposition}
The results analogous to Theorems \ref{u-basis} and \ref{u-alg}
for the complex projective and hyperbolic spaces seem to be of
interest and are still to be obtained.

\section{Main results: valuations on manifolds.}\label{main-results}
Let $X$ be a smooth manifold. Let $n=\dim X$. We assume also for
simplicity that $X$ is countable at infinity, i.e. $X$ can be
presented as a countable union of compact sets. Let us denote by
$\cp(X)$ the family of all compact submanifolds with corners (for
the notion of manifold with corners, see e.g. the book
\cite{melrose}, Chapter 1). We will be interested in finitely
additive measures on $\cp(X)$ which satisfy some additional
conditions of continuity (even some smoothness).
\begin{remark}
The class $\cp(X)$ is neither closed under finite unions nor under
finite intersections. Thus the notion of finite additivity should
be explained. Roughly, finite additivity holds whenever it makes
sense (see Section 2.2 of Part II \cite{part2} for the details).
\end{remark}
Let $P\in \cp(X)$ be a compact submanifold with corners. Let us
remind the definition of the characteristic cycle of $P$ denoted
by $CC(P)$. For any point $x\in P$ let $T_xP\subset T_xX$ denote
the tangent cone of $P$ at the point $x$. It is defined as
follows:
$$T_xP:=\{\xi\in T_xX| \mbox{ there exists a } C^1-\mbox{map }
\gamma\colon [0,1]\to P \mbox{ such that }\gamma (0)=x \mbox { and
} \gamma'(0)=\xi\}.$$ $T_xP$ is a convex polyhedral cone, and if
$P$ has no corners then $T_xP$ is the usual tangent space at $x$.
Let $(T_xP)^\circ \subset T^*_xX$ denote the dual cone. Recall
that the dual cone $C^o$ of a convex cone $C$ in a linear space
$W$ is defined by
$$C^o:=\{y\in W^*|\, y(x)\geq 0 \mbox{ for any } x\in C\}.$$
Define the characteristic cycle
\begin{eqnarray*}
CC(P):=\bigcup_{x\in P}(T_xP)^\circ.
\end{eqnarray*}
Then it is well known that $CC(P)$ has the following properties:

(1) $CC(P)\subset T^*X$ is an $n$-dimensional submanifold with
singularities;

(2) $CC(P)$ is Lagrangian, $\RR_{>0}$-invariant;

(3) if $X$ is oriented then $CC(P)$ is an $n$-cycle, i.e. $\pt
(CC(P))=0$.

\begin{remark}
(1) If $P$ has no corners then $CC(P)$ is the usual co-normal
bundle.

(2) The notion of characteristic cycle (or almost equivalent
notion of normal cycle) is well known. First the notion of normal
cycle was introduced by Wintgen \cite{wintgen}, and then studied
further by Z\"ahle \cite{zahle87} by the tools of geometric
measure theory. Characteristic cycles of subanalytic sets of real
analytic manifolds were introduced by J. Fu \cite{fu-94} using the
tools of geometric measure theory, and independently by Kashiwara
(see \cite{kashiwara-schapira}, Chapter 9) using the tools of
sheaf theory. J. Fu's article \cite{fu-94} develops an approach to
define the normal cycle for more general sets than subanalytic or
convex ones (see Theorem 3.2 in \cite{fu-94}).
\end{remark}

\begin{definition}
Let $\phi\colon \cp(X)\to \CC$ be a finitely additive measure. We
say that $\phi$ is {\itshape continuous} if for any uniformly
bounded sequence $\{P_N\}\subset \cp(X)$ and $P\in \cp(X)$ such
that $CC(P_N)\overset{flat}{\to} CC(P)$ one has
$$\phi(P_N)\to \phi(P).$$
\end{definition}
\begin{remark}
In the above definition the convergence is understood in the sense
of local flat convergence of currents on $T^*X$. For the
definition of this notion we refer to the book \cite{federer}.
Here we only notice that it is well known that if $\{K_N\}$ is a
uniformly bounded sequence of convex compact subsets in $\RR^n$
(or more generally, compact uniformly bounded subsets with reach
at least $\delta
>0$)  then $K_N\to K$ in the Hausdorff metric iff
$CC(K_N)\overset{flat}{\to} CC(K)$.
\end{remark}

The above assumption of continuity of measures is the most
important condition. It turns out mostly for technical reasons
that one should impose on valuations some other conditions. We
will skip here the precise definitions since most of these
conditions are technical and their necessity is not very clear for
the moment. Let us denote by $V^\infty(X)$ the set of smooth
measures which are called {\itshape smooth valuations}.
$V^\infty(X)$ is a linear space and has a natural nuclear
Fr\'echet topology. This space is the main object we are going to
discuss. For understanding of this survey it is possible to accept
the description of all smooth valuations given in Proposition
\ref{integ} below as a {\itshape definition} of smooth valuations.
\begin{examples}
(1) Any smooth density on $X$ belongs to $V^\infty(X)$.

(2) The Euler characteristic belongs to $V^\infty(X)$.
\end{examples}

Recall that we denote by $Val^\infty(\RR^n)$ the space of smooth
translation invariant convex valuations on $\RR^n$ in the sense of
Section \ref{overview}. Let us denote by $(V^\infty(\RR^n))^{tr}$
the subspace of $V^\infty(\RR^n)$ of translation invariant
valuations.
\begin{proposition}\label{comparison}
Restriction to convex subsets of $\RR^n$ defines the map
$$(V^\infty(\RR^n))^{tr}\to Val^\infty(\RR^n).$$
This map is an isomorphism.
\end{proposition}
\begin{remark}
In other words, this proposition says that the elements of
$Val^\infty(\RR^n)$ extend canonically to a class of not
necessarily convex sets.
\end{remark}
The next proposition provides a description of $V^\infty(X)$. We
will assume for simplicity that the manifold $X$ is oriented
though this is not strictly necessary, and the result can be
appropriately generalized to any smooth manifold.
\begin{proposition}\label{integ}
Let $X$ be an oriented manifold. Let $\phi\in V^\infty(X)$. There
exists a $C^\infty$-smooth differential $n$-form $\omega$ on
$T^*X$ which has the support compact relative to the canonical
projection $T^*X\to X$, and such that for any $P\in \cp(X)$
\begin{eqnarray}\label{inte}
\phi(P)=\int _{CC(P)}\ome.
\end{eqnarray}
And vice versa, any expression of the above form is a smooth
valuation.
\end{proposition}
\begin{remark}\label{rmk-kernel}
(1) For a given valuation $\phi$, the above form $\ome$ is highly
non-unique. For instance, for the Euler characteristic the above
form is well known and it was constructed by Chern \cite{chern}.
The construction depends on an extra choice of a Riemannian metric
on $X$.

(2) The forms defining the zero valuation were described in a
recent preprint by Bernig and Br\"ocker \cite{bernig-brocker} by
an explicitly written system of differential and integral
equations. In particular, they realized key role of the Rumin
operator on differential forms on contact manifolds \cite{rumin}
for this problem. We refer to \cite{bernig-brocker} for the
precise statements.

(3) The fact that the expressions of the form (\ref{inte}) are
smooth valuations heavily uses the tools from geometric measure
theory (see \cite{fu-89}-\cite{fu-94}). The converse statement
uses the Irreducibility Theorem (Theorem \ref{irrthm}) in
combination with the Casselman-Wallach theorem (see e.g.
\cite{wallach}) from representation theory.

\end{remark}
\begin{proposition}[Sheaf property]
The correspondence for any open subset $U\subset X$
$$U\mapsto V^\infty(U)$$ with the obvious restriction maps is a sheaf on $X$ which we denote by $\cv^\infty_X$.
\end{proposition}
\def\cvx{\cv^\infty_X}
Let us denote by $Val(TX)$ the (infinite dimensional) vector
bundle over $X$ such that its fiber over a point $x\in X$ is equal
to the space $Val^\infty(T_xX)$ of smooth translation invariant
convex valuations on $T_xX$. By McMullen's theorem it has a
grading by the degree of homogeneity:
$Val^\infty(TX)=\oplus_{i=0}^nVal^\infty_i(TX).$
\begin{theorem}
There exists a canonical filtration of $V^\infty(X)$ by closed
subspaces $$V^\infty(X)=W_0\supset W_1\supset\dots\supset W_n$$
such that the associated graded space $\oplus _{i=0}^n
W_i/W_{i+1}$ is canonically isomorphic to the space of smooth
sections $C^\infty(X,Val^\infty_i(TX))$.
\end{theorem}
\begin{remark}
(1) For $i=n$ the above isomorphism means that $W_n$ coincides
with the space of smooth densities on $X$.

(2) For $i=0$ the above isomorphism means that $W_0/W_1$ is
canonically isomorphic to the space of smooth functions
$C^\infty(X)$.

(3) Actually $U\mapsto W_i(U)$ defines a subsheaf $\cw_i$ of
$\cvx$.
\end{remark}
Next, $V^\infty(X)$ carries a very important and non-trivial
multiplicative structure (extending to the product on translation
invariant valuations discussed in Theorem \ref{poincare}).
\begin{theorem}
There exists a canonical product $V^\infty(X)\times V^\infty(X)\to
V^\infty(X)$ which is

(1) continuous;

(2) commutative and associative;

(3) the filtration $W_\bullet$ is compatible with it:
$$W_i\cdot W_j\subset W_{i+j};$$

(4) the Euler characteristic $\chi$ is the unit in the algebra
$V^\infty(X)$;

(5) this product commutes with restrictions to open and closed
submanifolds; in particular $\cvx$ is a sheaf of filtered
algebras.
\end{theorem}
In order to formulate one of the most interesting properties of
this product, let us observe that the space of smooth valuations
with compact support $V^\infty_c(X)$ admits a continuous
integration functional
\begin{eqnarray}\label{intfunc}
\int\colon V^\infty_c(X)\to \CC
\end{eqnarray}
given by $\phi\mapsto \phi(X)$. Consider the bilinear map
$$V^\infty(X)\times V^\infty_c(X)\to \CC$$
defined by $(\phi, \psi)\mapsto \int\phi\cdot \psi$.
\begin{theorem}\label{selfduality}
This bilinear form is a perfect pairing. In other words, the
induced map
$$V^\infty(X)\to (V^\infty_c(X))^*$$
is injective and has a dense image (with respect to the weak
topology in $(V^\infty_c(X))^*$).
\end{theorem}
The proof of this theorem uses the full statement of the
Irreducibility Theorem for translation invariant convex
valuations. We call Theorem \ref{selfduality} the Selfduality
Property of valuations.
\begin{definition}
Let us define $V^{-\infty}(X):=(V^\infty_c(X))^*$. Elements of
this space will be called {\itshape generalized valuations}.
\end{definition}
Thus $V^{-\infty}(X)$ can be considered as a completion of
$V^\infty(X)$ with respect to the weak topology. Roughly speaking,
one may say that the space of valuations is essentially self-dual
(at least when $X$ is compact). For any open subsets $U\subset
V\subset X$ one has the natural restriction map $V^{-\infty}(V)\to
V^{-\infty}(U)$ dual to the imbedding
$V^\infty_c(U)\hookrightarrow V^\infty_c(V)$. The assignment
$U\mapsto V^{-\infty}(U)$ is a sheaf. This sheaf is denoted by
$\cv^{-\infty}_X$.

\begin{theorem}
There exists a canonical automorphism of the algebra of smooth
valuations $$\sigma\colon V^\infty(X)\to V^\infty(X)$$ such that

(1) $\sigma^2=Id,$ i.e. $\sigma$ is involutive;

(2) $\sigma$ is continuous;

(3) $\sigma$ preserves the filtration $W_\bullet$, namely
$\sigma(W_i)=W_i$;

(4) for any smooth translation invariant valuation $\phi$ on
$\RR^n$ one has $$(\sigma \phi)(K)=(-1)^{\deg \phi}\phi(-K)$$
where $\deg \phi$ denotes the degree of homogeneity of $\phi$.

(5) The involution $\sigma$ extends (uniquely) by continuity to
$V^{-\infty}(X)$ in the weak topology.

(6) $\sigma$ commutes with restrictions to open subsets, thus it
defines the involution of sheaves $\cv_X^\infty$ and
$\cv^{-\infty}_X$.
\end{theorem}
\begin{remark}
We call $\sigma$ the Euler-Verdier involution.
\end{remark}

Let us discuss valuations on real analytic manifolds. On these
manifolds the space of constructible functions imbeds canonically
as a dense subspace of the space of generalized valuations, and it
is useful to compare the properties of the space of valuations
with the more familiar properties of the space of constructible
functions.

\begin{definition}[\cite{kashiwara-schapira}, \S9.7]\label{construct-function}
An integer valued function $f\colon
X\to \ZZ$ is called constructible if

1) for any $m\in \ZZ$ the set $f^{-1}(m)$ is subanalytic;

2) the family of sets $\{f^{-1}(m)\}_{m\in \ZZ}$ is locally
finite.
\end{definition}

Clearly the set of constructible $\ZZ$-valued functions is a ring
with pointwise multiplication. As in \cite{kashiwara-schapira} we
denote this ring by $CF(X)$. Define
\begin{eqnarray}\label{constr-def}
\cf:=CF(X)\otimes_\ZZ \CC.
\end{eqnarray}
Thus $\cf$ is a subalgebra of the $\CC$-algebra of complex valued
functions on $X$. In the rest of the article the elements of $\cf$
will be called {\itshape constructible functions}.
\def\cfx{\cf(X)}

 Let $\cf_c(X)$ denote the subspace of $\cf(X)$ of
{\itshape compactly supported} constructible functions. Clearly
$\cf_c(X)$ is a subalgebra of $\cf(X)$ (without unit if $X$ is
non-compact).

For a subset $P\subset X$ let us denote by $\One_P$ the indicator
function of $P$, namely
\begin{eqnarray*}
\One_P(x)=\left\{\begin{array}{ccc}
                       1&\mbox{ if }&x\in P\\
                       0&\mbox{ if }&x\not\in P.
                  \end{array}\right.
\end{eqnarray*}

For any constructible function $f$ one can define (see
\cite{kashiwara-schapira}, \cite{fu-94}) its characteristic cycle
$CC(f)$ so that for any subanalytic subset $P\subset X$ which is a
compact subanalytic submanifold with corners one has
$CC(\One_P)=CC(P)$.

We have a canonical linear map
\begin{eqnarray}\label{xi}
\Xi\colon \cf\to V^{-\infty}(X)
\end{eqnarray}
 which is uniquely characterized by
the property
$$<\Xi(\One_P),\phi>=\phi(P)$$
for any closed subanalytic subset $P$ and any $\phi\in
V^\infty_c(X)$ (strictly speaking, such a subset $P$ does not
belong to $\cp(X)$ in general, so this construction must be
justified; we refer to Section 8 of \cite{part4} for the details).
The map $\Xi$ turns out to be injective and has a dense image in
the weak topology. Thus we get the following imbeddings of dense
subspaces of $V^{-\infty}(X)$
\begin{eqnarray}\label{imbedd}
\cf\inj V^{-\infty}(X)\hookleftarrow V^\infty(X).
\end{eqnarray}

First one can extend the filtration $\{W_\bullet\}$ on
$V^\infty(X)$ to $V^{-\infty}(X)$ by taking closures in the weak
topology of each $W_i$. By \cite{part4} the restriction of this
filtration on $V^{-\infty}(X)$ back to $V^\infty(X)$ coincides
with the original filtration $\{W_\bullet\}$. Now let us restrict
the filtration obtained on $V^{-\infty}(X)$ to $\cf$ and denote
the induced filtration by $\{\tilde W_\bullet\}$. Then the
filtration $\{\tilde W_\bullet\}$ on $\cf$ coincides with the
filtration by the codimension of the support. More precisely we
have
\begin{proposition} Let $X$ be a real analytic manifold of dimension $n$.
For any $i=0,1,\dots,n$
$$\tilde W_i=\{f\in \cf|\, \codim (\supp f)\geq i\}$$
\end{proposition}

For compactly supported valuations and constructible functions we
have imbeddings analogous to (\ref{imbedd})
\begin{eqnarray}\label{imbeddc}
\cf_c\inj V^{-\infty}_c(X)\hookleftarrow V^\infty_c(X).
\end{eqnarray}
Then the integration functional $\int\colon V^\infty_c(X)\to\CC$
extends (uniquely) by continuity in the weak topology to a linear
functional
$$\int\colon V^{-\infty}_c(X)\to \CC.$$
\begin{proposition}
Let $X$ be a real analytic manifold. The restriction of the
integration functional from $V^{-\infty}_c(X)$ to $\cf_c$
coincides with the integration with respect to the Euler
characteristic.
\end{proposition}
Recall that the integration with respect to the Euler
characteristic is a linear functional $\cf_c\to \CC$ which is
uniquely characterized by the property $\One_P\mapsto \chi(P)$ for
any compact subanalytic subset $P\subset X$.

\begin{remark}\label{rmk-product}
Observe that $\cf$ has the natural structure of a $\CC$-algebra
with the usual pointwise multiplication. In some sense, this
product should correspond to the canonical product on smooth
valuations $V^\infty(X)$ discussed above. It would be interesting
to make this statement rigorous. Notice however that it seems to
be very unlikely that the product extends to the whole space
$V^{-\infty}(X)$ of generalized valuations.
\end{remark}
\begin{proposition}
The restriction of the Euler-Verdier involution $\sigma\colon
V^{-\infty}(X)\to V^{-\infty}(X)$  to $\cf$ coincides with the
classical Verdier involution times $(-1)^n$ where $n=\dim X$.
\end{proposition}
For the definition of the classical Verdier involution on the
constructible functions we refer to the book by Kashiwara-Schapira
\cite{kashiwara-schapira}, Chapter 9. Here we notice only that if
$P$ is a closed subanalytic submanifold with corners then
$\sigma(\One_P)=(-1)^{n-\dim P}\One_{\textrm{int}\, P}$ where
$\textrm{int}\, P$ denotes the relative interior of $P$.

\section{Appendix: a short guide to \cite{part1}, \cite{part2},
\cite{part3}, \cite{part4}.}\label{guide} Let us describe very
briefly the structure of the articles \cite{part1}, \cite{part2},
\cite{part3}, \cite{part4}.

Part I \cite{part1} still works only with {\itshape convex}
valuations on linear spaces. There one introduces a class of
smooth {\itshape convex} valuations playing a key role in
subsequent articles since it serves as a bridge between convex
valuations and general ones. Another main issue of Part I
\cite{part1} is a construction of a product on this class of
smooth convex valuations. It generalizes the previous construction
of the product \cite{alesker-gafa-04} for smooth convex valuations
which are in addition {\itshape polynomial with respect to
translations}. This generalization was based on the author's
earlier article \cite{alesker-int}.

Part II \cite{part2} introduces smooth valuations on general
smooth manifolds. A number of descriptions of this notion are
presented, and the comparison with smooth convex valuations
introduced in Part I \cite{part1}. In particular, the description
in terms of integration with respect to the normal
(characteristic) cycle is presented. This description uses a
number of results on convex valuations (including the
Irreducibility Theorem) in combination with results from geometric
measure theory (discussed in greater detail in Part III
\cite{part3}), and the Casselman-Wallach theorem from
representation theory. Also the Euler-Verdier involution on smooth
valuations was introduced in \cite{part2}.

The main goal of Part III \cite{part3} is to extend the product on
smooth convex valuations from Part I \cite{part1} to smooth
valuations on general manifolds. Roughly it works as follows.
Choosing a coordinate atlas for $X$, one uses the product of
valuations on $\RR^n$, defined by the construction of Part I
\cite{part1}, to define the product locally. Then one shows that
the products obtained on each coordinate patch coincide on
pairwise intersections, and that the result does not depend on the
choice of atlas. This step uses geometric measure theory. One
proves commutativity, associativity, and continuity with respect
to the natural topology of this product. Also in Part III
\cite{part3} one reviews and proves a number of relevant results
on normal cycles and geometric measure theory (following mostly
\cite{fu-89}-\cite{fu-94}) which were crucial in the description
of smooth valuations in terms of integration with respect to the
normal cycle discussed in Part II \cite{part2}.

The goal of Part IV \cite{part4} is twofold. First one studies
further the properties of the product of smooth valuations. In
particular, one proves that the filtration $\{W_\bullet\}$ on
smooth valuations is compatible with the product, namely $W_i\cdot
W_j\subset W_{i+j}$. Then one shows that the Euler-Verdier
involution is an automorphism of the algebra of smooth valuations.
Then one introduces the integration functional on compactly
supported valuations and proves the Selfduality Property (Theorem
\ref{selfduality} in this text). The second main point of Part IV
\cite{part4} is introducing the notion of generalized valuations
and establishing of basic properties of them.

\end{document}